\documentclass[10pt,final]{article}
\usepackage{a4}
\usepackage{amsmath}%
\usepackage{amstext}%
\usepackage{amssymb}%
\usepackage{amsthm}
\usepackage{comment}
\usepackage{listings}

\usepackage{url}

\usepackage[final]{graphicx}
\usepackage{subcaption}
\usepackage[framemethod=tikz]{mdframed}

\usepackage[margin=0.5in]{geometry}

\usepackage{graphicx}
\usepackage{epstopdf}

\usepackage{xcolor}

\newtheorem{theorem}{Theorem}

\newtheorem{rem}[theorem]{Remark}
\newenvironment{remark}{\rem\rm}{\endrem}

\newcounter{unnumber}

\newtheorem{prf}[unnumber]{Proof}
\renewenvironment{proof}{\prf\rm}{\hfill{$\blacksquare$}\endprf}
\newcommand{\R}{\mathbb{R}}%
\newcommand{\e}{\varepsilon}%
\newcommand{\ol}{\overline}%

\newcommand{\n}{{\nabla}}

\newcommand{\To}{\longrightarrow}
\def\a{\alpha}

\def\b{\beta}
\def\e{\epsilon}

\def\t{\theta}
\def\g{\gamma}

\def\<{\langle}
\def\>{\rangle}
\DeclareMathOperator*\id{Id}%

\DeclareMathOperator*\zer{zer}

\DeclareMathOperator*\cl{cl}
\author{Ern\"o Robert Csetnek\thanks{Faculty of Mathematics, University of Vienna, Oskar-Morgenstern-Platz 1, 1090 Vienna, Austria, email: robert.csetnek@univie.ac.at} \and Szil\'{a}rd Csaba L\'{a}szl\'{o} \thanks{Technical University of Cluj-Napoca, Department of Mathematics, Str. Memorandumului nr. 28, 400114 Cluj-Napoca, Romania, email: slaszlo@math.utcluj.ro.}}

\title{Strong convergence and fast rates for systems with Tikhonov regularization  \footnote{This work was supported by a grant of the Romanian Ministry of Research, Innovation and Digitization, CNCS - UEFISCDI, project number PN-III-P1-1.1-TE-2021-0138, within PNCDI III. }}

\begin{document}

\maketitle

\noindent \textbf{Abstract.}  We introduce and investigate the asymptotic behaviour of the trajectories of a second order dynamical system with Tikhonov regularization for solving a monotone equation with single valued, monotone and continuous operator acting on a real Hilbert space. We consider a vanishing damping which is correlated with the Tikhonov parameter and which is in line with recent developments in the literature on this topic. A correction term which involves the time derivative of the operator along the trajectory is also involved in the system and makes the link with Newton and Levenberg-Marquardt type methods. We obtain strong convergence of the trajectory to the minimal norm solution and fast convergence rates for the velocity and a quantity involving the operator along the trajectory. The rates are very closed to the known fast convergence results for systems without Tikhonov regularization, the novelty with respect to this is that we also obtain strong convergence of the trajectory to the minimal norm solution. As an application we introduce a primal-dual dynamical system for solving linearly constrained convex optimization problems, where the strong convergence of the trajectories is highlighted together with fast convergence rates for the feasibility measure and function values.

\vspace{1ex}

 \vspace{1ex}


\section{Introduction}
The aim of this paper is to approach the monotone equation
\begin{equation}\label{ME}\tag{ME}
Ax=0,
\end{equation}
by means of dynamical systems, where $A:\mathcal{H}\to\mathcal{H}$ is a monotone continuous operator and $\mathcal{H}$ is a real Hilbert space endowed with a scalar product $\<\cdot,\cdot\>$ and induced norm $\|\cdot\|$. We assume that the solution set $S:=\{x\in\mathcal{H}:Ax=0\}$ is non-empty.

To this purpose we associate to \eqref{ME} the second order dynamical system

 \begin{align}\label{DynSys1}\tag{DS}
\begin{cases}
\ddot{x}(t) + \frac{\a}{t^q} \dot{x}(t) + \b t^q\frac{d}{dt}A(x(t))+\left(1-\frac{\g}{t^s}\right)A(x(t))+\frac{c}{t^{2q+s}}x(t)=0,\\
x(t_0) = u_0\in\mathcal{H}, \,
\dot{x}(t_0) = v_0\in\mathcal{H},
\end{cases}
\end{align}
where, $\a>1,\,q,\,\b,\, \g,\, s,\,c>0$ with $0<q+s<1$ and  the starting time in \eqref{DynSys1} is $t_0 > 0.$ Let us notice the presence of the Tikhonov term $\frac{c}{t^{2q+s}}x(t)$ in (DS) which will induce strong convergence of the trajectory to the minimal norm solution and also fast rates, as described in the next section.

Let us make the connection to other dynamical systems from the literature for monotone inclusions. First, notice that in \cite{BCL-jma} a similar dynamical system with Tikhonov regularization has been investigated in case the operator $A$ maximally monotone. In this case the operator $A$ can be regularized through its Yosida regularization (which is a cocoercove operator) and combine this with the powerful techniques of Tikhonov regularization, see \cite{BCL-jma}. However, let us mention that there are many applications, for example for saddle point problems where the conditions characterizing these points lead to a monotone inclusion where the operator does not satisfy the cocoercivity property. One can consider also the rotation operator $(x,y)\to(-y,x)$ which is skew, Lipschitz, but not cocoercive. Moreover, a dynamical system where the operator $A$ is evaluated by forward steps (and not backward like in the case of Yosida regularization) is more tractable.

Further, note that \eqref{DynSys1} is strongly related to the dynamical system without Tikhonov regularization, as studied in \cite{BCN} in connection with the monotone equation \eqref{ME}. Indeed, in \cite{BCN} the following dynamical system was considered:
 \begin{align}\label{BCNDS}
\begin{cases}
\ddot{x}(t) + \frac{\a}{t} \dot{x}(t) + \b(t)\frac{d}{dt}A(x(t))+\frac12\left(\dot{\b}(t)+\frac{\a}{t}\b(t)\right)A(x(t)))=0,\\
x(t_0) = u_0\in\mathcal{H}, \,
\dot{x}(t_0) = v_0\in\mathcal{H},
\end{cases}
\end{align}
where $\a\ge 2$ and $\b$ is a continuously differentiable nondecreasing function which satisfies
$$0<\sup_{t\ge t_0}\frac{t\dot{\b}(t)}{\b(t)}\le\a-2.$$
The authors in \cite{BCN} showed that the trajectory $x(t)$ generated by \eqref{BCNDS} converges to a solution of \eqref{ME} in the weak topology of $\mathcal{H}.$ Further, the following estimates hold:
$$\|\dot{x}(t)\|=o\left(\frac{1}{t}\right)\mbox{ as }t\to+\infty\mbox{ and }\|A(x(t))\|=o\left(\frac{1}{t\b(t)}\right)\mbox{ as }t\to+\infty.$$

Let us remark that we cannot allow such a general coefficient $\b(t)$ in \eqref{DynSys1}. Indeed, this is due some technical reasons that appear in the proof of strong convergence of the trajectories generated by \eqref{DynSys1} to the minimum norm zero of $A$. That is also the reason why we use the damping parameter $\frac{\a}{t^q},\,q<1$ and we introduce the Tikhonov regularization term $\frac{c}{t^{2q+s}}x(t).$ Note that by introducing Tikhonov regularization terms in dynamical systems and numerical algorithms lead to strong convergence of the generated trajectories/iterates and allows an a-priori control to which solution these trajectories/iterates converge. Usually we are interested in the minimum norm solution, that is, the projection of the origin on the solution set, see \cite{ABCR,L-jde,AL-nemkoz,L-coap, BCL-jma,L-CNSN,KL-amop,L-Siopt}.

However, according to Theorem \ref{strconverg}, beside the strong convergence of the trajectories generated by \eqref{DynSys1} to the minimum norm zero of $A$ we are able to obtain the rates
$$\|\dot{x}(t))\|=\mathcal{O}(t^{s-1})+\mathcal{O}\left(t^\frac{-2q-s}{2}\right)\mbox{ as }t\to+\infty$$
and
$$\|A(x(t))\|=\mathcal{O}(t^{s-1-q})+\mathcal{O}\left(t^\frac{-4q-s}{2}\right)\mbox{ as }t\to+\infty.$$
Note that in case $s=\frac{2(1-q)}{3}$ one gets
$$\|\dot{x}(t))\|=\mathcal{O}\left(t^{-\frac{1+2q}{3}}\right)\mbox{ as }t\to+\infty\mbox{ and }\|A(x(t))\|=\mathcal{O}\left(t^{-\frac{1+5q}{3}}\right)\mbox{ as }t\to+\infty.$$
Now, if one takes $\b(t)=\b_0t^q,\,\b_0>0$ in \eqref{BCNDS}, then the rates obtained in \cite{BCN} become
$$\|\dot{x}(t)\|=o\left(t^{-1}\right)\mbox{ as }t\to+\infty\mbox{ and }\|A(x(t))\|=o\left(t^{-1-q}\right)\mbox{ as }t\to+\infty$$
and since for $q$ close to 1 one has $-\frac{1+2q}{3}\thicksim -1$ and $-\frac{1+5q}{3}\thicksim -1-q$ in this case the rates obtained in the present paper and \cite{BCN} are comparable.

Let us mention here a recent contribution \cite[Section 3]{bn-arx}, where a second order dynamical system with strong convergence to the minimal norm solution and convergence rates has been investigated. The approach is different in comparison to our paper. Indeed, in \cite{bn-arx} the authors start with a first order system with Tikhonov regularization and then the results are extrapolated to second order systems by means of time rescaling and anchoring procedure.

Now let us discuss the damping considered in (DS) (i.e. the coefficient of the velocity $\dot x$). We emphasize that our damping $\frac{\a}{t^q}, q<1$ is slightly different from the consecrated damping $\frac{\a}{t}$ considered first in \cite{SBC} in the context of a dynamical system that models Nesterov's accelerated gradient method and used in a many papers, see \cite{AAD1,ACPR, ACR,APR,ALP,AP, AL-nemkoz,ALSiopt,BCL,BCL-AA,AL}. Recall that the dynamical system studied in \cite{SBC} in connection to the minimization problem of a smooth convex function $f$ has the form
$$(SBC)\,\,\,\,\ddot{x}(t)+\frac{\a}{t}\dot{x}(t)+\n f(x(t))=0,\,x(t_0)=u_0\in\mathcal{H},\,\dot{x}(t_0)=v_0\in\mathcal{H},\,\a\ge 3.$$
Two related systems are the dynamical systems introduced in \cite{APR} and \cite{ALP}, that is
$$(APR)\,\,\,\,\ddot{x}(t)+\frac{\a}{t}\dot{x}(t)+\b\n^2 f(x(t))\dot{x}(t)+\n f(x(t))=0,\,x(t_0)=u_0\in\mathcal{H},\,\dot{x}(t_0)=v_0\in\mathcal{H},\,\a\ge 3,\,\b>0$$
and
$$(ALP)\,\,\,\,\ddot{x}(t)+\frac{\a}{t}\dot{x}(t)+\n f\left(x(t)+\left(\g+\frac{\b}{t}\right)\dot{x}(t)\right)=0,\,x(t_0)=u_0\in\mathcal{H},\,\dot{x}(t_0)=v_0\in\mathcal{H},\,\a\ge 3,\,\b\in\R,\,\g\ge 0.$$

The systems $(SBC),\,(APR)$ and $(ALP)$ share a similar behaviour, in case $\a>3$ there is weak convergence of the generated trajectories to a minimizer of $f$ and the rates $f(x(t))-\min f=\mathcal{O}(t^{-2})$ as $t\to+\infty$ and $\|\dot{x}(t)\|=\mathcal{O}(t^{-1})$ as $t\to+\infty$ hold.

If one introduces the Tikhonov regularization term $\frac{c}{t^p}x(t),\,c>0,\,p>2$ in the above systems, then according to \cite{ACR,BCL,ALSiopt} the strong convergence result $\liminf_{t\to+\infty}\|x(t)-x^*\|=0$, where $x^*$ is the minimum norm minimizer of $f$, holds. However, neither in \cite{ACR} nor in \cite{BCL} rates for the potential energy $f(x(t))-\min f$ and velocity $\|\dot{x}(t)\|$ are provided together with the strong convergence of the trajectories.

Let us mention that a positive result in this context, namely $\liminf_{t\to+\infty}\|x(t)-x^*\|=0$ and the rates $f(x(t))-\min f=\mathcal{O}(t^{-2})$ as $t\to+\infty$ and $\|\dot{x}(t)\|=\mathcal{O}(t^{-1})$ as $t\to+\infty$ are provided in \cite{AL-nemkoz} for the dynamical system
$$(AL)\,\,\,\,\ddot{x}(t)+\frac{\a}{t}\dot{x}(t)+\n f(x(t))+\frac{c}{t^2}x(t)=0,\,x(t_0)=u_0\in\mathcal{H},\,\dot{x}(t_0)=v_0\in\mathcal{H},\,\a\ge 3.$$

Moreover, in order to improve the liminf result and obtain a full strong convergence result along with fast convergence rates one can follow a series of papers of Attouch and his collaborators on this subject, where the idea is to involve the Tikhonov parameter also in the damping coefficient. Indeed, \cite{L-jde}   deals with  the dynamical system
\begin{align}\label{DynSys}
&\ddot{x}(t) + \frac{\a}{t^q} \dot{x}(t) +\nabla f\left(x(t)\right) +\frac{c}{t^p}x(t)=0,\,
x(t_0) = u_0, \,
\dot{x}(t_0) = v_0,
\end{align}
where $t_0 > 0$, $0<q<1$ and $(u_0,v_0) \in \mathcal{H} \times \mathcal{H}$. According to \cite{L-jde} for an appropriate setting of the parameters  one has the convergence rates
$f(x(t))-\min f=o\left(t^{-1-q}\right),\,\|\dot{x}(t)\|=o\left(t^{-\frac{1+q}{2}}\right)$ as $t\to+\infty$, and the strong convergence of the trajectories to the minimal norm solution. Note that in case $q$ is close to 1, these rates are comparable with the consecrated rates $\mathcal{O}(t^{-2})$ and $\mathcal{O}(t^{-1})$, respectively. Moreover, the damping $\frac{\a}{t^{q}},\,q<1$ was successfully used in \cite{BCL-jma} in the context of a dynamical system associated to a monotone inclusion, where beside the strong convergence of the generated trajectories fast convergence rates were provided. Another successful use of this damping parameter is provided in \cite{L-coap} in the context of a dynamical system with a Tikhonov regularization related to $(ALP)$.

The main contributions of our paper are:

\begin{itemize}
 \item We provide a setting with the following results: we show that the trajectory of (DS) converges strongly to the minimal norm solution of the equation (ME) and show fast convergence rates for $\|\dot x(t)\|$ and $\|A(x(t))\|$ (Theorem 1)
 \item In Section 3 we provide an application for linearly constrained convex optimization problems: we introduce a primal-dual dynamical system with strong convergence of the trajectories and fast convergence rates for the feasibility measure $\|B(x(t))-b\|$ and objective function $|f(x(t))-f(x^*)|$ and compare the results with the literature on this subject (Theorem 2).
\end{itemize}

\section{Convergence analysis}

We investigate in this paper the following dynamical system with Tikhonov regularization term

\begin{mdframed}
 \begin{equation}\label{DynSys11}\tag{DS}
\ddot{x}(t) + \frac{\a}{t^q} \dot{x}(t) + \b t^q\frac{d}{dt}A(x(t))+\left(1-\frac{\g}{t^s}\right)A(x(t))+\frac{c}{t^{2q+s}}x(t)=0,\\
\end{equation}
with initial conditions $x(t_0) = u_0\in\mathcal{H}, \,
\dot{x}(t_0) = v_0\in\mathcal{H}.$
\end{mdframed}
where, $\a>1,\,q,\,\b,\, \g,\, s,\,c>0$ with $0<q+s<1$ and  the starting time in \eqref{DynSys11} is $t_0 > 0.$
Further, we assume that $c<\frac{8\a(\a-1)\g}{\a^2\b^2+8(\a-1)\b}.$

We are interested in a strong global solution of (DS), which is a trajectory $x:[t_0,+\infty)\to\mathcal{H}$ with the property that both $x$ and $\dot x$ are absolutely continuous on every compact interval $[t_0,T]$, with $T>t_0$. We assume that $t\to A(x(t))$ is absolutely continuous on every compact interval $[t_0,T]$, which implies that $t\to A(x(t))$ is almost everywhere differentiable on $[t_0,+\infty)$. Notice that if we impose Lipschitz continuity conditions on the operator $A$, then one can show existence and uniqueness of a strong global solution of (DS) (see for example \cite{BCN}).

Now let us recall some basics from the Tikhonov regularization theory which will be used in the investigations below, see \cite{BCL-jma}. Consider $t \mapsto\epsilon(t)$ a positive function such that $\e(t)\To 0$ as $t\To+\infty$. Let $x_{\varepsilon(t)}$ be the unique zero of the strongly monotone operator $A+\e(t)\id,$ that is $$Ax_{\e(t)} +\e(t)x_{\e(t)}= 0.$$ Then $x_{\e(t)}$ converges (as $\,t\to+\infty$) in the strong topology to $x^*$, which is the projection of $0$ onto the set $S:= \{x \in \mathcal{H} :  Ax=0\}$, i.e. the minimal norm solution. Furthermore, we have $\|x_{\e(t)}\|\leq\|x^*\|$.

Moreover, if $\epsilon:[t_0,+\infty)\to (0,+\infty)$ is a nonincreasing function of class $C^1$ then $t\to x_{\e(t)}$ is almost everywhere differentiable and we have the estimate
\begin{equation}\label{dx_t_epsilon}\left\|\frac{d}{dt} x_{\e_t}\right\|\le\left( -\frac{\dot{\e}(t)}{\e(t)}\right)\|x_{\e_t}\|\mbox{ for almost all } t\ge t_0.\end{equation}

In the following we will use a particular regularization, namely $\e(t)=\frac{c}{t^p}$, with $c,p>0$. In this case \eqref{dx_t_epsilon} becomes
\begin{equation}\label{dx_t}\left\|\frac{d}{dt} x_{t}\right\|\le\frac{p}{t}\|x_{t}\|\mbox{ for almost all } t\ge t_0,\end{equation}
where, by a slight abuse of notation, we denote $x_{t}$ as the unique  zero of the strongly monotone operator $A_{t}x:=Ax+\frac{c}{t^p}x,$ that is $$A(x_{t}) +\frac{c}{t^p}x_{t}= 0.$$

The main result concerning the properties of the trajectories generated by \eqref{DynSys11} is stated in the following theorem.

\begin{theorem}\label{strconverg} For some starting points $u_0,v_0\in\mathcal{H}$ let $x(t)$ the trajectory generated by \eqref{DynSys11}. Let $x^*$ be the minimum norm zero of the operator $A$ and $x_t$ the unique  zero of the strongly monotone operator $A_{t}x:=Ax+\frac{c}{t^{2q+s}}x.$ Then the following statements hold.
\begin{itemize}
\item[(i)]
$\|x(t)-x_t\|=\mathcal{O}(t^{q+s-1})+\mathcal{O}\left(t^\frac{-s}{2}\right),$
consequently, $x(t)$ converges strongly to $x^*$, that is,
$$\lim_{t\to+\infty}x(t)=x^*.$$
\end{itemize}
Further, we obtain the following rates:
\begin{itemize}
\item[(ii)]$\|A(x(t))\|=\mathcal{O}(t^{-1-q+s})+\mathcal{O}\left(t^\frac{-4q-s}{2}\right)$ as $t\to+\infty.$
\item[(iii)]$\|\dot{x}(t)\|=\mathcal{O}(t^{-1+s})+\mathcal{O}\left(t^\frac{-2q-s}{2}\right)$ as $t\to+\infty.$
\end{itemize}
\end{theorem}
\begin{proof}
First consider $b>0$ such that $2\a>b$ and for simplicity denote $\a(t)=\frac{\a}{t^q}$,\,$\b(t)=\b t^q$ and $\g(t)=1-\frac{\g}{t^s}$ and $\e(t)=\frac{c}{t^{2q+s}}.$ Let $c(t)=t^q$, $$v(t)=\frac12\|b(x(t)-x_t)+c(t)(2\dot{x}(t)+\b(t)A(x(t))\|^2$$ and
$$u(t)=bc(t)\b(t)\<A(x(t)),x(t)-x_t\>+\frac{bc(t)\b(t)\e(t)}{2}(\|x(t)\|^2-\|x(t)-x_t\|^2-\|x_t\|^2).$$
Note that $u(t)=bc(t)\b(t)\<A_{t}(x(t))-A_{t}(x_t),x(t)-x_t\>-bc(t)\b(t)\e(t)\|x(t)-x_t\|^2\ge0.$

We introduce the energy functional
\begin{align*}{E}(t)={v}(t)&+{u}(t)+\frac{c^2(t)\b^2(t)}{2}\|A(x(t))\|^2+2c^2(t)\e(t)\|x(t)\|^2\\
&+\frac{b(2\a(t)c(t)-b-2\dot{c}(t)+c(t)\b(t)\e(t))}{2}\|x(t)-x_t\|^2.
\end{align*}
Note that for $t$ big enough $b(2\a(t)c(t)-b-2\dot{c}(t)+c(t)\b(t)\e(t))\ge 0$, hence also $E(t)\ge 0.$

In the first part of the proof, which is very technical, we will show that \eqref{at} holds.

By making use of the system \eqref{DynSys11}, we get
\begin{align}\label{v1}
\dot{v}(t)&=(b+2\dot{c}(t)-2\a(t)c(t))\left<\dot{x}(t), b(x(t)-x_t)+c(t)(2\dot{x}(t)+\b(t)A(x(t))\right>\\
\nonumber&+(\dot{c}(t)\b(t)+c(t)\dot{\b}(t)-2c(t)\g(t))\left\<A(x(t)), b(x(t)-x_t)+c(t)(2\dot{x}(t)+\b(t)A(x(t))\right>\\
\nonumber&-c(t)\b(t)\left\<\frac{d}{dt}A(x(t)), b(x(t)-x_t)+c(t)(2\dot{x}(t)+\b(t)A(x(t))\right>\\
\nonumber&-\left\<2c(t)\e(t)x(t)+b\frac{d}{dt} x_t, b(x(t)-x_t)+c(t)(2\dot{x}(t)+\b(t)A(x(t))\right>.
\end{align}

Further we have

\begin{align}\label{v11}
\left<\dot{x}(t), b(x(t)-x_t)+c(t)(2\dot{x}(t)+\b(t)A(x(t))\right>=2c(t)\|\dot{x}(t)\|^2+c(t)\b(t)\<A(x(t)),\dot{x}(t)\>+b\<\dot{x}(t),x(t)-x_t\>,
\end{align}
and

\begin{align}\label{v12}
\left\<A(x(t)), b(x(t)-x_t)+c(t)(2\dot{x}(t)+\b(t)A(x(t))\right>&=c(t)\b(t)\|A(x(t))\|^2+2c(t)\<A(x(t)),\dot{x}(t)\>\\
\nonumber&+b\<A(x(t)), x(t)-x_t\>
\end{align}
By using the fact that due to the monotonicity of $A$ one has $\<\frac{d}{dt}A(x(t)), \dot{x}(t)\>\ge 0$ we derive
\begin{align}\label{v13}
&-c(t)\b(t)\left\<\frac{d}{dt}A(x(t)), b(x(t)-x_t)+c(t)(2\dot{x}(t)+\b(t)A(x(t))\right>\le -bc(t)\b(t)\left\<\frac{d}{dt}A(x(t)), x(t)-x_t\right>\\
\nonumber&+\frac{(c^2(t)\b^2(t))'}{2}\|A(x(t))\|^2-\frac{d}{dt}\left(\frac{c^2(t)\b^2(t)}{2}\|A(x(t))\|^2\right).
\end{align}

Further,
\begin{align}\label{v14}
&-\left\<2c(t)\e(t)x(t)+b\frac{d}{dt} x_t, b(x(t)-x_t)+c(t)(2\dot{x}(t)+\b(t)A(x(t))\right>=bc(t)\e(t)(\|x_t\|^2-\|x(t)\|^2-\|x(t)-x_t\|^2\\
\nonumber&-4c^2(t)\e(t)\<\dot{x}(t),x(t)\>-2c^2(t)\e(t)\b(t)\<A(x(t)),x(t)\>-b\left\<\frac{d}{dt} x_t, b(x(t)-x_t)+c(t)(2\dot{x}(t)+\b(t)A(x(t))\right>.
\end{align}

Combining \eqref{v1}, \eqref{v11}, \eqref{v12}, \eqref{v13} and \eqref{v14}  we get
\begin{align}\label{v2}
\dot{v}(t)&+\frac{d}{dt}\left(\frac{c^2(t)\b^2(t)}{2}\|A(x(t))\|^2\right)\le 2(b+2\dot{c}(t)-2\a(t)c(t))c(t)\|\dot{x}(t)\|^2\\
\nonumber&+\left((b+2\dot{c}(t)-2\a(t)c(t))c(t)\b(t)+2(\dot{c}(t)\b(t)+c(t)\dot{\b}(t)-2c(t)\g(t))c(t)\right)\<A(x(t)),\dot{x}(t)\>\\
\nonumber&+(2\dot{c}(t)\b(t)+2c(t)\dot{\b}(t)-2c(t)\g(t))c(t)\b(t)\|A(x(t))\|^2\\
\nonumber&-bc(t)\e(t)\|x(t)-x_t\|^2-bc(t)\e(t) \|x(t)\|^2+bc(t)\e(t)\|x_t\|^2\\
\nonumber&+b(\dot{c}(t)\b(t)+c(t)\dot{\b}(t)-2c(t)\g(t))\<A(x(t)),x(t)-x_t\>-2c^2(t)\b(t)\e(t)\<A(x(t)),x(t)\>\\
\nonumber&+b(b+2\dot{c}(t)-2\a(t)c(t))\<\dot{x}(t),x(t)-x_t\>-4c^2(t)\e(t)\<\dot{x}(t),x(t)\>\\
\nonumber&-bc(t)\b(t)\left\<\frac{d}{dt}A(x(t)),x(t)-x_t\right\>-b\left\<\frac{d}{dt} x_t, b(x(t)-x_t)+c(t)(2\dot{x}(t)+\b(t)A(x(t))\right>.
\end{align}

Let us  give some further estimates:
\begin{align*}
b(b+2\dot{c}(t)-2\a(t)c(t))\<\dot{x}(t),x(t)-x_t\>&=\frac{d}{dt}\left(\frac{b(b+2\dot{c}(t)-2\a(t)c(t))}{2}\|x(t)-x_t\|^2\right)\\
&+b(b+2\dot{c}(t)-2\a(t)c(t))\left\<\frac{d}{dt}{x_t},x(t)-x_t\right\>-b\ddot{c}(t)\|x(t)-x_t\|^2,
\end{align*}

\begin{align*}
-4c^2(t)\e(t)\<\dot{x}(t),x(t)\>&=-\frac{d}{dt}(2c^2(t)\e(t)\|x(t)\|^2)+(2c^2(t)\e(t))'\|x(t)\|^2
\end{align*}
and
\begin{align*}
-2c^2(t)\b(t)\e(t)\<A(x(t)),x(t)\>&=-2c^2(t)\b(t)\e(t)\<A(x(t)),x(t)-x_t\>-2c^2(t)\b(t)\e(t)\<A(x(t)),x_t\>.
\end{align*}

Hence, \eqref{v2} becomes
\begin{align}\label{v3}
\dot{v}(t)&+\frac{d}{dt}\left(\frac{c^2(t)\b^2(t)}{2}\|A(x(t))\|^2\right)+\frac{d}{dt}(2c^2(t)\e(t)\|x(t)\|^2)+\frac{d}{dt}\left(\frac{b(2\a(t)c(t)-b-2\dot{c}(t)}{2}\|x(t)-x_t\|^2\right)\le\\
\nonumber&2(b+2\dot{c}(t)-2\a(t)c(t))c(t)\|\dot{x}(t)\|^2\\
\nonumber&+\left((b+2\dot{c}(t)-2\a(t)c(t))c(t)\b(t)+2(\dot{c}(t)\b(t)+c(t)\dot{\b}(t)-2c(t)\g(t))c(t)\right)\<A(x(t)),\dot{x}(t)\>\\
\nonumber&+(2\dot{c}(t)\b(t)+2c(t)\dot{\b}(t)-2c(t)\g(t))c(t)\b(t)\|A(x(t))\|^2\\
\nonumber&-(bc(t)\e(t)+b\ddot{c}(t))\|x(t)-x_t\|^2+((2c^2(t)\e(t))'-bc(t)\e(t)) \|x(t)\|^2+bc(t)\e(t)\|x_t\|^2\\
\nonumber&+\left(b(\dot{c}(t)\b(t)+c(t)\dot{\b}(t)-2c(t)\g(t))-2c^2(t)\b(t)\e(t)\right)\<A(x(t)),x(t)-x_t\>-2c^2(t)\b(t)\e(t)\<A(x(t)),x_t\>\\
\nonumber&+b(2\dot{c}(t)-2\a(t)c(t))\left\<\frac{d}{dt} x_t,x(t)-x_t\right\>-2bc(t)\left\<\frac{d}{dt} x_t, \dot{x}(t)\right>-bc(t)\b(t)\left\<\frac{d}{dt} x_t,A(x(t))\right>\\
\nonumber&-bc(t)\b(t)\left\<\frac{d}{dt}A(x(t)),x(t)-x_t\right\>.
\end{align}

Considering $u(t)$, one has
\begin{align}\label{u1}
\dot{u}(t)&=b(\dot{c}(t)\b(t)+c(t)\dot{\b}(t))\<A(x(t)),x(t)-x_t\>+bc(t)\b(t)\left\<A(x(t)),\dot{x}(t)-\frac{d}{dt}x_t\right\>\\
\nonumber&+bc(t)\b(t)\left\<\frac{d}{dt}A(x(t)),x(t)-x_t\right\>+\frac{d}{dt}\left(\frac{bc(t)\b(t)\e(t)}{2}\|x(t)\|^2\right)
-\frac{d}{dt}\left(\frac{bc(t)\b(t)\e(t)}{2}\|x(t)-x_t\|^2\right)\\
\nonumber&-\frac{d}{dt}\left(\frac{bc(t)\b(t)\e(t)}{2}\|x_t\|^2\right).
\end{align}

Now,
$$\frac{d}{dt}\left(\frac{bc(t)\b(t)\e(t)}{2}\|x(t)\|^2\right)=\frac{(bc(t)\b(t)\e(t))'}{2}\|x(t)\|^2+bc(t)\b(t)\e(t)\<\dot{x}(t),x(t)\>,$$
consequently \eqref{v3} and \eqref{u1} lead to

\begin{align}\label{uv1}
\dot{E}(t)&=\dot{v}(t)+\dot{u}(t)+\frac{d}{dt}\left(\frac{c^2(t)\b^2(t)}{2}\|A(x(t))\|^2\right)+\frac{d}{dt}\left(2c^2(t)\e(t)\|x(t)\|^2\right)
\\
\nonumber&+\frac{d}{dt}\left(\frac{b(2\a(t)c(t)-b-2\dot{c}(t)+c(t)\b(t)\e(t))}{2}\|x(t)-x_t\|^2\right)\le2(b+2\dot{c}(t)-2\a(t)c(t))c(t)\|\dot{x}(t)\|^2\\
\nonumber&+\left((2b+2\dot{c}(t)-2\a(t)c(t))c(t)\b(t)+2(\dot{c}(t)\b(t)+c(t)\dot{\b}(t)-2c(t)\g(t))c(t)\right)\<A(x(t)),\dot{x}(t)\>\\
\nonumber&+(2\dot{c}(t)\b(t)+2c(t)\dot{\b}(t)-2c(t)\g(t))c(t)\b(t)\|A(x(t))\|^2\\
\nonumber&-(bc(t)\e(t)+b\ddot{c}(t))\|x(t)-x_t\|^2+\left((2c^2(t)\e(t))'+\frac{(bc(t)\b(t)\e(t))'}{2}-bc(t)\e(t)\right) \|x(t)\|^2+bc(t)\e(t)\|x_t\|^2\\
\nonumber&+\left(b(2\dot{c}(t)\b(t)+2c(t)\dot{\b}(t)-2c(t)\g(t))-2c^2(t)\b(t)\e(t)\right)\<A(x(t)),x(t)-x_t\>-2c^2(t)\b(t)\e(t)\<A(x(t)),x_t\>\\
\nonumber&+b(2\dot{c}(t)-2\a(t)c(t))\left\<\frac{d}{dt} x_t,x(t)-x_t\right\>-2bc(t)\left\<\frac{d}{dt} x_t, \dot{x}(t)\right>-2bc(t)\b(t)\left\<\frac{d}{dt} x_t,A(x(t))\right>\\
\nonumber&+bc(t)\b(t)\e(t)\<\dot{x}(t),x(t)\>-\frac{d}{dt}\left(\frac{bc(t)\b(t)\e(t)}{2}\|x_t\|^2\right).
\end{align}

By using $\left\|\frac{d}{dt} x_t\right\|\le -\frac{\dot{\e}(t)}{\e(t)}\|x_t\|=\frac{q+2s}{t}\|x_t\|\mbox{ for almost every } t\ge t_0,$
we have the following estimates for every $s_1,s_2,s_3.s_4>0$.
\begin{align}\label{est1}
-2c^2(t)\b(t)\e(t)\<A(x(t)),x_t\>\le c^2(t)\b(t)\e(t)\left(\frac{s_1}{\e(t)}\|A(x(t))\|^2+\frac{\e(t)}{s_1}\|x_t\|^2\right),
\end{align}

\begin{align}\label{est2}
b(2\dot{c}(t)-2\a(t)c(t))\left\<\frac{d}{dt} x_t,x(t)-x_t\right\>\le b(\a(t)c(t)-\dot{c}(t))\left(s_2c(t)\e(t)\|x(t)-x_t\|^2+\frac{\dot{\e}^2(t)}{s_2c(t)\e^3(t)}\|x_t\|^2\right),
\end{align}

\begin{align}\label{est3}
-2bc(t)\left\<\frac{d}{dt} x_t, \dot{x}(t)\right>\le bc(t)\left(s_3\|\dot{x}(t)\|^2+\frac{\dot{\e}^2(t)}{s_3\e^2(t)}\|x_t\|^2\right)
\end{align}
and
\begin{align}\label{est4}
2bc(t)\b(t)\left\<\frac{d}{dt} x_t,A(x(t))\right>\le bc(t)\b(t)\left(s_4c(t)\|A(x(t))\|^2+\frac{\dot{\e}^2(t)}{s_4c(t)\e^2(t)}\|x_t\|^2\right).
\end{align}

Further, for every $s_5>0$ one has
\begin{align}\label{est5+}
&bc(t)\b(t)\e(t)\<\dot{x}(t),x(t)\>\le s_5bc(t)\|\dot{x}(t)\|^2+\frac{bc(t)\b^2(t)\e^2(t)}{4s_5}\|x(t)\|^2.
\end{align}
Now, by using the monotonicity of $A$ one obtains
\begin{align*}\<A(x(t)),x(t)-x_t\>&=\<A(x(t))-A(x_t),x(t)-x_t\>-\e(t)\<x_t,x(t)-x_t\>\ge-\e(t)\<x_t,x(t)-x_t\>\\
&=\frac{\e(t)}{2}(\|x_t\|^2+\|x(t)-x_t\|^2-\|x(t)\|^2).
\end{align*}
By taking into account that $2\dot{c}(t)\b(t)+2c(t)\dot{\b}(t)-2c(t)\g(t)\le 0$ for $t$ big enough we get

\begin{align}\label{est5}
&\left(b(2\dot{c}(t)\b(t)+2c(t)\dot{\b}(t)-2c(t)\g(t))-2c^2(t)\b(t)\e(t)\right)\<A(x(t)),x(t)-x_t\>\le\\
\nonumber&\left(b(\dot{c}(t)\b(t)+c(t)\dot{\b}(t)-c(t)\g(t))-c^2(t)\b(t)\e(t)\right)\e(t)\|x_t\|^2\\
\nonumber&+\left(b(\dot{c}(t)\b(t)+c(t)\dot{\b}(t)-c(t)\g(t))-c^2(t)\b(t)\e(t)\right)\e(t)\|x(t)-x_t\|^2\\
\nonumber&-\left(b(\dot{c}(t)\b(t)+c(t)\dot{\b}(t)-c(t)\g(t))-c^2(t)\b(t)\e(t)\right)\e(t)\|x(t)\|^2.
\end{align}

Hence, we have for almost every $t$ (where $t$ is large enough) that
\begin{align}\label{uv2}
\dot{E}(t)&\le (2b+s_3b+s_5b+4\dot{c}(t)-4\a(t)c(t))c(t)\|\dot{x}(t)\|^2\\
\nonumber&+\left((2b+2\dot{c}(t)-2\a(t)c(t))c(t)\b(t)+2(\dot{c}(t)\b(t)+c(t)\dot{\b}(t)-2c(t)\g(t))c(t)\right)\<A(x(t)),\dot{x}(t)\>\\
\nonumber&+(2\dot{c}(t)\b(t)+2c(t)\dot{\b}(t)-2c(t)\g(t)+s_1c(t)+bs_4c(t))c(t)\b(t)\|A(x(t))\|^2\\
\nonumber&+\left(b((\a-\dot{c}(t))s_2-1)c(t)\e(t)-b\ddot{c}(t)+\left(b((c(t){\b}(t))'-c(t)\g(t))-c^2(t)\b(t)\e(t)\right)\e(t)\right)\|x(t)-x_t\|^2\\
\nonumber&+\left((2c^2(t)\e(t))'+\frac{(bc(t)\b(t)\e(t))'}{2}+\frac{bc(t)\b^2(t)\e^2(t)}{4s_5}-bc(t)\e(t)-\left(b((c(t){\b}(t))'-c(t)\g(t))
-c^2(t)\b(t)\e(t)\right)\e(t)\right) \|x(t)\|^2\\
\nonumber&+\left(bc(t)\e(t)+\frac{c^2(t)\b(t)\e^2(t)}{s_1}+\frac{b(\a-\dot{c}(t))\dot{\e}^2(t)}{s_2c(t)\e^3(t)}
+\frac{bc(t)\dot{\e}^2(t)}{s_3\e^2(t)}+\frac{b\b(t)\dot{\e}^2(t)}{s_4\e^2(t)}\right)\|x_t\|^2\\
\nonumber&+\left(b((c(t){\b}(t))'-c(t)\g(t))-c^2(t)\b(t)\e(t)\right)\e(t)\|x_t\|^2-\frac{d}{dt}\left(\frac{bc(t)\b(t)\e(t)}{2}\|x_t\|^2\right).
\end{align}

Consider now $K>0$ and  denote $r=q+s$. Then,
\begin{align*}
\frac{K}{t^r}v(t)&\le \frac{Kb^2}{t^r}\|x(t)-x_t\|^2+ \frac{8Kc^2(t)}{t^r}\|\dot{x}(t)\|^2+\frac{2Kc^2(t)\b^2(t)}{t^r}\|A(x(t))\|^2
\end{align*}
and for every $s_6>0$ one has
\begin{align*}
\frac{K}{t^r}u(t)&=\frac{bKc(t)\b(t)}{t^r}\<A(x(t)),x(t)-x_t\>+\frac{bKc(t)\b(t)\e(t)}{2t^r}(\|x(t)\|^2-\|x(t)-x_t\|^2-\|x_t\|^2)\\
&\le \frac{s_6bKc^2(t)\b(t)}{2}\|A(x(t))\|^2+\frac{bK\b(t)}{2s_6t^{2r}}\|x(t)-x_t\|^2+\frac{bKc(t)\b(t)\e(t)}{2t^r}(\|x(t)\|^2-\|x(t)-x_t\|^2-\|x_t\|^2).
\end{align*}
Consequently,
\begin{align}\label{KrE}
\frac{K}{t^r}E(t)&\le\frac{8Kc^2(t)}{t^r}\|\dot{x}(t)\|^2+K\left(\frac{s_6bt^r+5\b(t)}{2t^r}\right)c^2(t)\b(t)\|A(x(t))\|^2+
K\left(\frac{2c^2(t)\e(t)}{t^r}+\frac{bc(t)\b(t)\e(t)}{2t^r}\right)\|x(t)\|^2\\
\nonumber&+K\left(\frac{b(2\a(t)c(t)+b-2\dot{c}(t))}{2t^r}+\frac{b\b(t)}{2s_6t^{2r}}\right)\|x(t)-x_t\|^2-\frac{bKc(t)\b(t)\e(t)}{2t^r}\|x_t\|^2.
\end{align}

Now, by adding \eqref{uv2} and \eqref{KrE} we get
\begin{align}\label{at}
\dot{E}(t)&+\frac{K}{t^r}E(t)\le R_1(t)\|\dot{x}(t)\|^2+R_2(t)\<A(x(t)),\dot{x}(t)\>+R_3(t)\|A(x(t))\|^2+R_4(t)\|x(t)-x_t\|^2+R_5(t) \|x(t)\|^2\\
\nonumber&+\left(bc(t)\e(t)+\frac{c^2(t)\b(t)\e^2(t)}{s_1}+\frac{b(\a-\dot{c}(t))\dot{\e}^2(t)}{s_2c(t)\e^3(t)}
+\frac{bc(t)\dot{\e}^2(t)}{s_3\e^2(t)}+\frac{b\b(t)\dot{\e}^2(t)}{s_4\e^2(t)}-\frac{bKc(t)\b(t)\e(t)}{2t^r}\right)\|x_t\|^2\\
\nonumber&+\left(b((c(t){\b}(t))'-c(t)\g(t))-c^2(t)\b(t)\e(t)\right)\e(t)\|x_t\|^2-\frac{d}{dt}\left(\frac{bc(t)\b(t)\e(t)}{2}\|x_t\|^2\right),
\end{align}
where the coefficients $R_1-R_5$ will be analyzed in what follows.  The aim is to show that \eqref{at1} below holds.

We use $c(t)=t^q,\,\b(t)=\b t^q,\,\g(t)=1-\frac{\g}{t^s},\,\b,\g,s>0,\,0<q+s= r<\min(1,q+2s),$ $1<\a$ and we have

\begin{align}\label{R1}
R_1(t)&=\left(2b+s_3b+s_5b+4\dot{c}(t)-4\a(t)c(t)+\frac{8Kc(t)}{t^r}\right))c(t)\\
\nonumber&=\left(2b+s_3b+s_5b+4qt^{q-1}-4\a+\frac{8Kt^q}{t^r}\right)t^q=(2b+s_3b+s_5b-4\a)t^q+\mathcal{O}(t^{q-s})
\end{align}

\begin{align}\label{R2}
R_2(t)&=\left((2b+2\dot{c}(t)-2\a(t)c(t))c(t)\b(t)+2(\dot{c}(t)\b(t)+c(t)\dot{\b}(t)-2c(t)\g(t))c(t)\right)\\
\nonumber&=(2b+2qt^{q-1}-2\a)\b t^{2q}+2(2q\b t^{2q-1}-2+2\g t^{q-s})t^q=((2b-2\a)\b-4)t^{2q}+6q\b t^{3q-1}+4\g t^{2q-s}
\end{align}

\begin{align}\label{R3}
R_3(t)&=(2\dot{c}(t)\b(t)+2c(t)\dot{\b}(t)-2c(t)\g(t)+s_1c(t)+bs_4c(t))c(t)\b(t)+K\left(\frac{s_6bt^r+5\b(t)}{2t^r}\right)c^2(t)\b(t)\\
\nonumber&=\left(2qt^{q-1}\b t^q+2q\b t^{2q-1}-2 t^q+\frac{2\g t^q}{t^s}+s_1t^q+bs_4t^q\right)\b t^{2q}+K\left(\frac{s_6bt^r+5\b t^q}{2t^r}\right)\b t^{3q}\\
\nonumber&=\left(-2+s_1+bs_4+\frac{Ks_6b}{2}\right)\b t^{3q}+\mathcal{O}(t^{3q-s})
\end{align}
Our next step is to show that for appropriately chosen $s_1,s_3,s_4,s_5,s_6$  we have
$$R_1(t)\|\dot{x}(t)\|^2+R_2(t)\<A(x(t)),\dot{x}(t)\>+R_3(t)\|A(x(t))\|^2<0,\mbox{ for }t\mbox{ big enough}.$$
 To this purpose is enough to show that $R_1(t),\,R_3(t)<0$ and $R_2^2(t)-4 R_1(t)R_3(t)<0$ for $t$ big enough.

Obviously, we can take $s_1,s_3,s_4,s_5,s_6$ as small as we need hence $R_3(t)<0$ for $t$ big enough and from $b<2\a$ we find that $R_1(t)<0$ for $t$ big enough. Further, the dominant coefficient, i.e. the coefficient of $t^{4q}$, in the expression $R_2^2(t)-4 R_1(t)R_3(t)$ is
\begin{align}\label{forneg}
&((2b-2\a)\b-4)^2-4(2b+s_3b+s_5b-4\a)\left(-2+s_1+bs_4+\frac{Ks_6b}{2}\right)=((2b-2\a)\b-4)^2+8(2b-4\a)\\
\nonumber&-4(2b-4\a)\left(-s_1+bs_4+\frac{Ks_6b}{2}\right)-4(s_3+s_5)b\left(-2+s_1+bs_4+\frac{Ks_6b}{2}\right).
\end{align}
 Now, if for $s_1,s_3,s_4,s_6\to 0$ we have $((2b-2\a)\b-4)^2+8(2b-4\a)+8s_5b<0$ for some $s_5>0$ then obviously we can fix $s_1,s_3,s_4,s_5,s_6>0$ such  that \eqref{forneg} is negative.

Let us show that one can have $((2b-2\a)\b-4)^2+8(2b-4\a)+8s_5b<0$.
Since $\a>1$ we take $b=\a$, hence $((2b-2\a)\b-4)^2+8(2b-4\a)+8s_5b=16(1-\a)+8s_5\a<0$ for all $0<s_5<\frac{2(\a-1)}{\a}.$
 Now, we fix $s_1,s_3,s_4,s_5,s_6$ such that
 $((2b-2\a)\b-4)^2-4(2b+s_3b+s_5b-4\a)\left(-2+s_1+bs_4+\frac{Ks_6b}{2}\right)<0$, and we get
$$R_1(t)\|\dot{x}(t)\|^2+R_2(t)\<A(x(t)),\dot{x}(t)\>+R_3(t)\|A(x(t))\|^2<0,\mbox{ for }t\mbox{ big enough}.$$


Further for $s_2=\frac{1}{\a}$ we have
\begin{align}\label{R4}
R_4(t)&=\left(b((\a-\dot{c}(t))s_2-1)c(t)\e(t)-b\ddot{c}(t)+\left(b((c(t){\b}(t))'-c(t)\g(t))-c^2(t)\b(t)\e(t)\right)\e(t)\right)\\
\nonumber&+K\left(\frac{b(2\a(t)c(t)+b-2\dot{c}(t))}{2t^r}+\frac{b\b(t)}{2s_6t^{2r}}\right)\\
\nonumber&=\left(\a((\a-qt^{q-1})s_2-1)t^q\e(t)+\a q(1-q)t^{q-2}+\left(\a (2\b qt^{2q-1}-t^q+\g t^{q-s})-\b t^{3q}\e(t)\right)\e(t)\right)\\
\nonumber&+K\left(\frac{\a(3\a-2q t^{q-1})}{2t^r}+\frac{\a\b t^q}{2s_6t^{2r}}\right)=\a(\a s_2-2)t^q\e(t)+\frac{3K\a^2}{2t^r}-\b t^{3q}\e^2(t)+\mathcal{O}(t^{q-s}\e(t))+\mathcal{O}(t^{q-2r})\\
\nonumber&=\left(-\a c+\frac{3K\a^2}{2}\right)t^{-q-s}+\mathcal{O}(t^{-q-2s}),
\end{align}
and  we conclude that $R_4(t)\le 0$ for the appropriate  $K<\frac{2 c}{3\a}$ and $t$ big enough.

Finally, we get

\begin{align}\label{R5}
R_5(t)&=\left((2c^2(t)\e(t))'+\frac{(bc(t)\b(t)\e(t))'}{2}+\frac{bc(t)\b^2(t)\e^2(t)}{4s_5}-bc(t)\e(t)-\left(b((c(t){\b}(t))'-c(t)\g(t))
-c^2(t)\b(t)\e(t)\right)\e(t)\right)\\
\nonumber&+K\left(\frac{2c^2(t)\e(t)}{t^r}+\frac{bc(t)\b(t)\e(t)}{2t^r}\right)
\\
\nonumber&=\frac{4+\a\b}{2}(t^{2q}\e(t))'+\frac{\a\b^2}{4s_5}t^{3q}\e^2(t)-\a t^q\e(t)-\left(\a(2\b qt^{2q-1}-t^q+\g t^{q-s})-\b t^{3q}\e(t)\right)\e(t))\\
\nonumber&+K\frac{4+\a\b}{2}t^{2q-r}\e(t)\\
\nonumber&=-\a\g t^{q-s}\e(t)+\frac{4+\a\b}{2}t^{2q}\dot{\e}(t)+\left(\frac{\a\b^2}{4s_5}+\b\right) t^{3q}\e^2(t)+K\frac{4+\a\b}{2}t^{2q-r}\e(t)+\mathcal{O}(t^{2q-1}\e(t))\\
\nonumber&= \left(-\a\g c+\left(\frac{\a\b^2}{4s_5}+\b\right)c^2+Kc\frac{4+\a\b}{2}\right)t^{-q-2s}+\mathcal{O}(t^{-s-1}).
\end{align}

We show first that one can choose $0<s_5<\frac{2(\a-1)}{\a}$ such that $-\a\g +\left(\frac{\a\b^2}{4s_5}+\b\right) c\le -\t $ for some $\t>0$ and $t$ big enough. Indeed, since $c<\frac{8\a(\a-1)\g}{\a^2\b^2+8(\a-1)\b}$ we consider $\t>0$ such that $c<\frac{8(\a-1)(\a\g-\t)}{\a^2\b^2+8(\a-1)\b}$. Then, obviously $c<\frac{\a\g-\t}{\b}$ hence one can take $0<\frac{\a\b^2}{4\left(\frac{(\a\g-\t)}{c}-\b\right)}<s_5<\frac{2(\a-1)}{\a}$. Consequently,
$\left(\frac{\a\b^2}{4s_5}+\b\right)<\frac{\a\g-\t}{c}$ and the claim follows. Now, just simply take $K<\min \left(\frac{2\t}{4+\a\b}, \frac{2 c}{3\a}\right)$ and  get that  $R_5(t)\le 0$ for $t$ big enough.

Let us deal now with the coefficient of $\|x_t\|^2$ in \eqref{at}. We use the fact that  $\|x_t\|$ is bounded since $\|x_t\|\le\|x^*\|$ and we get

\begin{align}\label{R6}
R_6(t)&=\left(bc(t)\e(t)+\frac{c^2(t)\b(t)\e^2(t)}{s_1}+\frac{b(\a-\dot{c}(t))\dot{\e}^2(t)}{s_2c(t)\e^3(t)}
+\frac{bc(t)\dot{\e}^2(t)}{s_3\e^2(t)}+\frac{b\b(t)\dot{\e}^2(t)}{s_4\e^2(t)}-\frac{bKc(t)\b(t)\e(t)}{2t^r}\right)\|x_t\|^2\\
\nonumber&+\left(b((c(t){\b}(t))'-c(t)\g(t))-c^2(t)\b(t)\e(t)\right)\e(t)\|x_t\|^2-\frac{d}{dt}\left(\frac{bc(t)\b(t)\e(t)}{2}\|x_t\|^2\right)\\
\nonumber&=\left(\a t^q\e(t)+\frac{\b}{s_1} t^{3q}\e^2(t)+\frac{\a (\a-qt^{q-1})\dot{\e}^2(t)}{s_2t^q\e^3(t)}
+\frac{\a t^q\dot{\e}^2(t)}{s_3\e^2(t)}+\frac{\a\b t^q\dot{\e}^2(t)}{s_4\e^2(t)}-\frac{\a K\b t^{2q}\e(t)}{2t^r}\right)\|x_t\|^2\\
\nonumber&+\left(\a(2\b qt^{2q-1}-t^q+\g t^{q-s})-\b t^{3q}\e(t)\right)\e(t)\|x_t\|^2-\frac{d}{dt}\left(\frac{\a\b t^{2q}\e(t)}{2}\|x_t\|^2\right)\\
\nonumber&=\left(\b\left(\frac{1}{s_1}-1\right)t^{3q}\e^2(t)+ \frac{\a(\a-qt^{q-1})\dot{\e}^2(t)}{s_2t^q\e^3(t)}
+\frac{\a t^q\dot{\e}^2(t)}{s_3\e^2(t)}+\frac{\a \b t^q\dot{\e}^2(t)}{s_4\e^2(t)}-\frac{\a K\b t^{2q}\e(t)}{2t^r}\right)\|x_t\|^2\\
\nonumber&+\left(\a(2\b qt^{2q-1}+\g t^{q-s})\right)\e(t)\|x_t\|^2-\frac{d}{dt}\left(\frac{\a\b t^{2q}\e(t)}{2}\|x_t\|^2\right)\\
\nonumber&\le C_1 t^{r-2}+C_2 t^{q-2}-C_3t^{q-2r}+C_4t^{-q-2s}-\frac{d}{dt}\left(\frac{b\b t^{2q}\e(t)}{2}\|x_t\|^2\right)\\
\nonumber&\le C_5 t^{q+s-2}+C_6 t^{-q-2s}-\frac{d}{dt}\left(\frac{b\b t^{-s}}{2}\|x_t\|^2\right),
\end{align}
for some constants $C_1-C_6>0$ and $t$ big enough.

Consequently, \eqref{at} yields that there exists $t_1\ge t_0$ such that
\begin{align}\label{at1}
\dot{E}(t)&+\frac{K}{t^{q+s}}E(t)\le C_5 t^{q+s-2}+C_6 t^{-q-2s}-\frac{d}{dt}\left(\frac{b\b t^{-s}}{2}\|x_t\|^2\right),
\end{align}
for almost every $t\ge t_1.$

Now, by multiplying \eqref{at1} with $e^{\frac{K}{1-q-s}t^{1-q-s}}$  and using the fact that for $v,W>0,\,u\in\R$ there exists $C>0$ such that
$$t^{u}e^{Wt^v}\le C\frac{d}{dt}\left(t^{u-v+1}e^{Wt^v}\right)$$
we get
\begin{equation}\label{th}
\frac{d}{dt}\left(E(t)e^{\frac{K}{1-q-s}t^{1-q-s}}\right)\le C \frac{d}{dt}\left(t^{2q+2s-2}e^{\frac{K}{1-q-s}t^{1-q-s}}\right)+C \frac{d}{dt}\left(t^{-s}e^{\frac{K}{1-q-s}t^{1-q-s}}\right)-e^{\frac{K}{1-q-s}t^{1-q-s}}\frac{d}{dt}\left(\frac{b\b t^{-s}}{2}\|x_t\|^2\right),
\end{equation}
for some $C>0$ and for almost every $t\ge t_1.$

Now, by integrating \eqref{th} on an interval $[t_1,T],\,T>t_1$ we get for some $C_0,C_0',C',K'>0$ that

\begin{align}\label{th1}
E(T)&\le C T^{2q+2s-2}+C T^{-s}-\frac{b\b T^{-s}}{2}\|x_T\|^2+e^{-\frac{K}{1-q-s}tT^{1-q-s}}\left(\int_{t_1}^T K' t^{-2s-q}\|x^*\|^2  e^{\frac{K}{1-q-s}t^{1-q-s}} dt+C'\right)\\
\nonumber&\le C_0T^{2q+2s-2}+C_0 T^{-s}.
\end{align}

Now taking into account the form of $E(t)$ we get at once that
$$\|x(t)-x_t\|=\mathcal{O}(t^{q+s-1})+\mathcal{O}\left(t^\frac{-s}{2}\right),\mbox{ as }t\to+\infty$$
consequently, we have
$$\lim_{t\to+\infty}x(t)=x^*.$$
Further
$$\|A(x(t))\|=\mathcal{O}(t^{-1-q+s})+\mathcal{O}\left(t^\frac{-4q-s}{2}\right),\mbox{ as }t\to+\infty.$$

Since $v(t)\le C_0't^{2q+2s-2}+C_0' t^{-s}$  for $t$ big enough, we obtain that
$\frac{v(t)}{t^{\max(2q+2s-2,-s)}}\le C''_0$ for some constant $C''_0.$ Now taking into account the form of $v(t)=\frac12\|b(x(t)-x_t)+c(t)(2\dot{x}(t)+\b(t)A(x(t))\|^2$ and the previously stated rates for $\|x(t)-x_t\|$ and $\|A(x(t))\|$ we obtain

$$\|\dot{x}(t))\|=\mathcal{O}(t^{s-1})+\mathcal{O}\left(t^\frac{-2q-s}{2}\right),\mbox{ as }t\to+\infty.$$
\end{proof}

\section{Application to optimization problems with linear constraints}

The aim of this section is to present a direct consequence of the results obtained for linearly constrained convex optimization problems.

Consider the optimization problem
$$(P)\,\,\,\,\min_{x\in\mathcal{H}} f(x)\mbox{ s.t. }Bx=b,$$
where $f:\mathcal{H}\to\R$ is a smooth convex function and $B:\mathcal{H}\to\mathcal{G}$ is a continuous linear operator ($\mathcal{H}, \mathcal{G}$ are real Hilbert spaces). Further we assume that
$b\in sri B(\mathcal{H}).$ Here the strong relative interior of a convex set $C\subseteq \mathcal{G}$ is  $sri C =\{u\in\mathcal{G}: \mbox{cone}(C-u) = \cl(\mbox{span}(C-u))\}.$

According to \cite{BC}, $x^*\in\mathcal{H}$ is a solution of $(P)$ if and only if
\begin{equation}\label{lancar}
Bx^*=b\mbox{ and there exists }y^*\in\mathcal{G}\mbox{ such that }B^*y^*=-\n f(x^*),
\end{equation}
where $B^*:\mathcal{G}\to\mathcal{H}$ is the adjoint operator.
If we introduce the lagrangian
$l:\mathcal{H}\times\mathcal{G}\to\R,\,\,l(x,y)=f(x)+\<y,Bx-b\>=f(x)+\<B^* y,x\>-\<y,b\> $, then $(x^*,y^*)$ is a saddle point of the lagrangian, that is
$$l(x^*,y)\le l(x^*,y^*)\le l(x,y^*),\,\forall (x,y)\in\mathcal{H}\times \mathcal{G}.$$
Note that $y^*$ is a solution of the dual problem
$$(D)\,\,\,\max_{y\in\mathcal{G}}\{-f^*(-B^*y)-\<y,b\>\}.$$

Hence, $(x^*,y^*)$ satisfy  the Fermat rule, which reads as
$$\n_x l(x,y)=\n f(x)+B^* y=0,\,\n_y l(x,y)=Bx-b=0.$$
 In other words $(x^*,y^*)$ satisfies \eqref{lancar}.

 Consequently, if we consider the operator $A:\mathcal{H}\times\mathcal{G}\to \mathcal{H}\times\mathcal{G},\,\, A(x,y)=(\n_x l(x,y), -\n_y l(x,y))=(\nabla f(x)+B^*y,b-Bx)$, then the above mentioned Fermat rule is exactly $\eqref{ME}.$ Notice that since $l$ is convex in $x$ and concave in $y$, the operator $A$ is maximally monotone.

 In the view of the above mentioned facts, consider the problem $(P)$ and associate to $(P)$ the following primal-dual second order dynamical system:

 \begin{mdframed}
 \begin{align}\label{primaldual}
 \begin{cases}
\ddot{x}(t) + \frac{\a}{t^q} \dot{x}(t) + \b t^q\Big(\n^2 f(x(t))\dot x(t)+B^*(\dot y(t))\Big)+\left(1-\frac{\g}{t^s}\right)(\n f(x(t))+B^*y(t))+\frac{c}{t^{2q+s}}x(t)=0,\\
\ddot{y}(t) + \frac{\a}{t^q} \dot{y}(t) - \b t^qB(\dot x(t))+\left(1-\frac{\g}{t^s}\right)(b-Bx(t))+\frac{c}{t^{2q+s}}y(t)=0,\\
(x(t_0),y(t_0)) = (u_0^1,u_0^2)\in\mathcal{H}\times\mathcal{G}, \,
(\dot{x}(t_0),\dot{y}(t_0)) = (v_0^1,v_0^2)\in\mathcal{H}\times\mathcal{G}.
\end{cases}
 \end{align}
  \end{mdframed}
where we assume that $f$ is $C^2$. Moreover, we impose $\a>1,\,q,\,\b,\, \g,\, s,\,c>0$ with $0<q+s<1$ and the starting time in \eqref{DynSys11} is $t_0 > 0.$
Further, we assume that $c<\frac{8\a(\a-1)\g}{\a^2\b^2+8(\a-1)\b}.$

Note that if $(x^*,y^*)$ is the minimum norm zero of the operator $A(x,y)=(\n f(x)+B^* y, b-Bx)$, then $x^*$ is a solution of $(P)$ and $y^*$ is a solution of $(D)$ and $\|(x^*,y^*)\|\le\|(\ol x,\ol y)\|$ for every $(\ol x,\ol y)\in \zer A.$

Applying Theorem \ref{strconverg} to the dynamical system \eqref{primaldual} we get:

\begin{theorem}\label{pdstrconverg} For some starting points $(u_0^1,u_0^2),(v_0^1,v_0^2)\in\mathcal{H}\times\mathcal{G}$ let $(x(t),y(t))$ the trajectory generated by \eqref{primaldual}. Let $(x^*,y^*)$ be the minimum norm zero of the operator $A(x,y)=(\n f(x)+B^* y, b-Bx)$. Then the following statements hold.
\begin{itemize}
\item[(i)]
$\lim_{t\to+\infty}x(t)=x^*\mbox{ and }\lim_{t\to+\infty}y(t)=y^*$, i.e. strong convergence of the trajectories.
\end{itemize}
Further, we obtain the following rates:
\begin{itemize}
\item[(ii)]$\|\n f(x(t))+B^*y(t)\|,\|Bx(t)-Bx^*\|\in\mathcal{O}(t^{-1-q+s})+\mathcal{O}\left(t^\frac{-4q-s}{2}\right)$ as $t\to+\infty.$
\item[(iii)]$\|\dot{x}(t)\|,\|\dot{y}(t)\|\in\mathcal{O}(t^{-1+s})+\mathcal{O}\left(t^\frac{-2q-s}{2}\right)$ as $t\to+\infty.$
\item[(iv)] $|f(x(t))-f(x^*)|\in\mathcal{O}(t^{-1-q+s})+\mathcal{O}\left(t^\frac{-4q-s}{2}\right)$ as $t\to+\infty.$
\end{itemize}
\end{theorem}

\begin{proof} The statements in (i)-(iii) are a direct consequence of Theorem \ref{strconverg}.

(iv) We make use of the fact that the trajectories $x,y$ are bounded and combine this with the gradient inequality for convex functions. We have:
\begin{align*} f(x(t))-f(x^*)\leq & \ \< \nabla f(x(t)),x(t)-x^* \> \\
 = & \ \< \nabla f(x(t))+B^*y(t),x(t)-x^* \> - \< B^*y(t),x(t)-x^*\>\\ \leq & \ M_1 \|\nabla f(x(t))+B^*y(t)\| - \< y(t),Bx(t)-Bx^*\>\\
 \leq & \ M_1 \|\nabla f(x(t))+B^*y(t)\| + M_2\|Bx(t)-Bx^*\|,
\end{align*}
where $M_1,M_2>0$ and
\begin{align*} f(x(t))-f(x^*)\geq &  \ \< \nabla f(x^*),x(t)-x^* \>\\
= & \ \< -B^*y^*,x(t)-x^* \> \\
= & \ -\< y^*,Bx(t)-Bx^* \> \\
\geq  & \ -\|y^*\|\cdot \|B(x(t))-Bx^*\|.
\end{align*}
The conclusion follows from (ii).
\end{proof}

\begin{remark} Notice that we obtain strong convergence of the trajectories as well as rates for the feasibility measure $\|B(x(t))-b\|$ and function value $|f(x(t))-f(x^*)|$. If we make a similar choice as in the introduction, namely $s=\frac{2(1-q)}{3}$, then the aforementioned rates are of order $\mathcal{O}\left(t^{-\frac{1+5q}{3}}\right)\mbox{ as }t\to+\infty$, with $q<1$ arbitrarily close to $1$. Hence the rates are of order $\mathcal{O}\left(\frac{1}{t^{2-\delta}}\right)\mbox{ as }t\to+\infty$, where $\delta>0$ can be arbitrarily chosen. This means that we are very closed to the standard accelerated rates obtained in \cite{bn} for linearly constrained optimization problems with the advantage that we obtain strong convergence of the trajectories in the framework of Tikhonov regularization.

\end{remark}




\begin{thebibliography}{999}




\bibitem{AAD1}{V. Apidopoulos, J.-F. Aujol,  Ch. Dossal},
{\it The differential inclusion modeling the FISTA algorithm and optimality of convergence rate in the case $b \leq 3$}, SIAM J. Optim.  28(1)  (2018),  551---574.

\bibitem{ACPR} { H. Attouch, Z. Chbani, J. Peypouquet, P. Redont}, \textit{Fast convergence of inertial dynamics and algorithms with asymptotic vanishing viscosity}, Mathematical Programming 168 (1-2) (2018),  123--175.


\bibitem{ACR} {H. Attouch,  Z. Chbani, H. Riahi},
{ \it Combining fast inertial dynamics for convex optimization with Tikhonov regularization},
J. Math. Anal. Appl 457 (2018),  1065--1094.


\bibitem{AL-nemkoz} H. Attouch, S. L\'aszl\'o, {\it Convex optimization via inertial algorithms with vanishing Tikhonov regularization: fast convergence to the minimum norm solution}, Math Meth Oper Res 99, 307–347 (2024)

\bibitem{ALSiopt} C.D. Alecsa, S.C. L\'aszl\'o, {\it Tikhonov Regularization of a Perturbed Heavy Ball System with Vanishing Damping}, SIAM Journal on Optimization, 31(4), 2921-2954, 2021

\bibitem{ALP} C.D. Alecsa, S.C. L\' aszl\' o, T. Pin\c ta, {\it An Extension of the Second Order Dynamical System that Models Nesterov’s Convex Gradient Method}, Applied Mathematics and Optimization, 84, 1687-1716, 2021



\bibitem{AL} H. Attouch, S.C. László, {\it Continuous Newton-like Inertial Dynamics for Monotone Inclusions}, Set-Valued Var. Anal 29, 555–581 (2021)

\bibitem{AP} H. Attouch, J.  Peypouquet, {\it Convergence of inertial dynamics and proximal algorithms governed by maximal monotone operators},  Mathematical Programming   {174}(1-2), 391-432 (2019)

\bibitem{APR} H. Attouch, J. Peypouquet, P. Redont, {\it Fast convex optimization via inertial dynamics with Hessian driven damping}, J. Differ. Equ. 261(10), 5734–5783 (2016)


\bibitem{ABCR} {H. Attouch, A. Balhag, Z. Chbani, H. Riahi},  {\it Damped inertial dynamics with vanishing Tikhonov regularization: Strong asymptotic convergence towards the minimum norm solution}, Journal of Differential Equations  311, 29-58 (2022)

\bibitem{BC} Bauschke, H., Combettes, P.L.: {Convex Analysis and Monotone Operator Theory in Hilbert spaces}. CMS Books in Mathematics. Springer   (2011)


\bibitem{BCL} R.I. Bo\c t, E.R. Csetnek, S.C. L\'aszl\'o,     \textit{Tikhonov regularization of a second order dynamical system with Hessian damping}, Math. Program. 189  (2021), 151–186.

\bibitem{BCL-AA} R.I. Bo\c t,  E.R. Csetnek,  S.C. L\' aszl\' o, {\it A second order dynamical approach with variable damping to nonconvex smooth minimization}, Applicable Analysis (2018). https://doi.org/10.1080/00036811.2018.1495330

\bibitem{BCL-jma} { R.I. Bo\c t, E.R. Csetnek, S.C. L\'aszl\'o}, {\it On the strong convergence of continuous Newton-like inertial dynamics with Tikhonov regularization for monotone inclusions}, Journal of Mathematical Analysis and Applications, 530(2),  2024

\bibitem{BCN} { R.I. Bo\c t, E.R. Csetnek, D.K. Nguyen}, {\it Fast Optimistic Gradient Descent Ascent (OGDA) Method in Continuous and Discrete Time}, Found Comput Math 2023, https://doi.org/10.1007/s10208-023-09636-5

\bibitem{bn} {R. I. Bo\c t, D.-K. Nguyen}, {\it Improved convergence rates and trajectory convergence for primal-dual dynamical systems with vanishing damping}, Journal of Differential Equations 303, 369-406, 2021

\bibitem{bn-arx} {R. I. Bo\c t, D.-K. Nguyen}, {\it Tikhonov regularization of monotone operator flows not only ensures strong convergence of the trajectories but also speeds up the vanishing of the residuals}, https://arxiv.org/abs/2406.00852, (2024)

\bibitem{KL-amop} { M. Karapetyants, S.C. L\'aszl\'o}, {\it  A Nesterov type algorithm with double Tikhonov regularization: fast convergence of the function values and strong convergence to the minimal norm solution}, Appl Math Optim 90, 17 (2024)



\bibitem{L-jde} {S.C. L\'aszl\'o}, {\it On the strong convergence of the trajectories of a Tikhonov regularized second order dynamical system with asymptotically vanishing damping}, Journal of Differential Equations 362, 355-381 (2023)

\bibitem{L-coap} { S.C. L\'aszl\'o}, {\it Solving convex optimization problems via a second order dynamical system with implicit Hessian damping and Tikhonov regularization}, Comput Optim Appl (2024). https://doi.org/10.1007/s10589-024-00620-5

\bibitem{L-CNSN} S.C. L\'aszl\'o, {\it On the convergence of an inertial proximal algorithm with a Tikhonov regularization term}, https://arxiv.org/abs/2302.02115  (2023)

\bibitem{L-Siopt} { S.C. L\'aszl\'o}, {\it A proximal-gradient inertial algorithm with Tikhonov regularization: strong convergence to the minimal norm solution}, arXiv:2407.10350 (2024)

 \bibitem{SBC}{\sc W.  Su,  S. Boyd,  E. J. Cand\`es}, { \it
A Differential Equation for Modeling Nesterov's
Accelerated Gradient Method: Theory and Insights},  Journal of Machine Learning Research 17(153) (2016), 1-43.
NIPS, December 2014.

\end{thebibliography}
\end{document}